\documentclass{article}%
\usepackage{amssymb}
\usepackage{amsmath}
\usepackage{amsfonts}
\usepackage{graphicx}%
\providecommand{\U}[1]{\protect\rule{.1in}{.1in}}
\newtheorem{theorem}{Theorem}[section]

\newtheorem{corollary}[theorem]{Corollary}

\newtheorem{proposition}[theorem]{Proposition}
\newtheorem{remark}[theorem]{Remark}

\newenvironment{proof}[1][Proof]{\noindent\textbf{#1.} }{\ \rule{0.5em}{0.5em}}
\begin{document}

\author{Vadim E. Levit\\Ariel University Center of Samaria, Israel\\levitv@ariel.ac.il
\and Eugen Mandrescu\\Holon Institute of Technology, Israel\\eugen\_m@hit.ac.il}
\title{Critical independent sets and K\"{o}nig--Egerv\'{a}ry graphs}
\date{}
\maketitle

\begin{abstract}
A set $S$ of vertices is \textit{independent} in a graph $G$, and we write
$S\in\mathrm{Ind}(G)$, if no two vertices from $S$ are adjacent, and
$\alpha(G)$ is the cardinality of an independent set of maximum size, while
$\mathrm{core}(G)$ denotes the intersection of all maximum independent sets
\cite{levm3}.

$G$ is called a \textit{K\"{o}nig--Egerv\'{a}ry graph} if its order equals
$\alpha(G)+\mu(G)$, where $\mu(G)$ denotes the size of a maximum matching. The
number $def(G)=\left\vert V\left(  G\right)  \right\vert -2\mu(G)$ is the
\textit{deficiency} of $G$ \cite{lovpl}.

The number $d(G)=\max\{\left\vert S\right\vert -\left\vert N(S)\right\vert
:S\in\mathrm{Ind}(G)\}$ is the \textit{critical difference} of $G$. An
independent set $A$ is \textit{critical} if $\left\vert A\right\vert
-\left\vert N(A)\right\vert =d(G)$, where $N(S)$ is the neighborhood of $S$,
and $\alpha_{c}(G)$ denotes the maximum size of a critical independent set
\cite{Zhang}.

In \cite{Larson2009} it was shown that $G$ is K\"{o}nig--Egerv\'{a}ry graph if
and only if there exists a maximum independent set that is also critical,
i.e., $\alpha_{c}(G)=\alpha(G)$.

In this paper we prove that:

\emph{(i)} $d(G)=\left\vert \mathrm{core}(G)\right\vert -\left\vert
N(\mathrm{core}(G))\right\vert =\alpha(G)-\mu(G)=def\left(  G\right)  $ hold
for every K\"{o}nig--Egerv\'{a}ry graph $G$;

\emph{(ii)} $G$ is K\"{o}nig--Egerv\'{a}ry graph if and only if each maximum
independent set of $G$ is critical.


\end{abstract}

\section{Introduction}

Throughout this paper $G=(V,E)$ is a finite, undirected, loopless and without
multiple edges graph with vertex set $V=V(G)$ and edge set $E=E(G)$. If
$X\subset V$, then $G[X]$ is the subgraph of $G$ spanned by $X$. By $G-W$ we
mean the subgraph $G[V-W]$ , if $W\subset V(G)$. For $F\subset E(G)$, by $G-F$
we denote the partial subgraph of $G$ obtained by deleting the edges of $F$,
and we use $G-e$, if $W$ $=\{e\}$. The neighborhood of a vertex $v\in V$ is
the set $N(v)=\{w:w\in V$ \ \textit{and} $vw\in E\}$, while $N(A)=\cup
\{N(v):v\in A\}$ and $N[A]=A\cup N(A)$ for $A\subset V$.

A set $S\subseteq V(G)$ is \textit{independent} if no two vertices from $S$
are adjacent, and by $\mathrm{Ind}(G)$ we mean the set of all the independent
sets of $G$. An independent set of maximum size will be referred to as a
\textit{maximum independent set} of $G$, and the \textit{independence number
}of $G$ is $\alpha(G)=\max\{\left\vert S\right\vert :S\in\mathrm{Ind}(G)\}$.

Let us denote the set $\{S:S$ \textit{is a maximum independent set of} $G\}$
by $\Omega(G)$, and let \textrm{core}$(G)=\cap\{S:S\in\Omega(G)\}$
\cite{levm3}. A set $A\subseteq V(G)$ is a \textit{local maximum independent
set} of $G$ if $A\in\Omega(G[N[A]])$ \cite{LevMan2}.

\begin{theorem}
\cite{NemhTro}\label{th1} Every local maximum independent set of a graph is a
subset of a maximum independent set.
\end{theorem}

A matching (i.e., a set of non-incident edges of $G$) of maximum cardinality
$\mu(G)$ is a \textit{maximum matching}, and a \textit{perfect matching} is
one covering all vertices of $G$.

It is well-known that
\[
\lfloor\left\vert V\right\vert /2\rfloor+1\leq\alpha(G)+\mu(G)\leq\left\vert
V\right\vert
\]
hold for any graph $G=(V,E)$. If $\alpha(G)+\mu(G)=\left\vert V\right\vert $,
then $G$ is called a \textit{K\"{o}nig-Egerv\'{a}ry graph}. We attribute this
definition to Deming \cite{dem}, and Sterboul \cite{ster}. These graphs were
studied in \cite{bourpull,korach,levm2,levm4,LevMan3,lov,lovpl,pulleybl}, and generalized
in \cite{bourhams1,pasdema}.

\begin{figure}[h]
\setlength{\unitlength}{1cm}\begin{picture}(5,1.3)\thicklines
\multiput(1.5,0)(1,0){3}{\circle*{0.29}}
\put(2.5,1){\circle*{0.29}}
\put(1.5,0){\line(1,0){2}}
\put(2.5,0){\line(0,1){1}}
\put(3.5,0){\line(-1,1){1}}
\put(1,0.5){\makebox(0,0){$H_{1}$}}
\multiput(5,0)(1,0){4}{\circle*{0.29}}
\put(7,1){\circle*{0.29}}
\put(6,1){\circle*{0.29}}
\put(8,1){\circle*{0.29}}
\put(5,0){\line(1,0){3}}
\put(6,1){\line(1,0){1}}
\put(7,0){\line(0,1){1}}
\put(5,0){\line(1,1){1}}
\put(7,0){\line(1,1){1}}
\put(4.5,0.5){\makebox(0,0){$H_{2}$}}
\multiput(10,0)(1,0){3}{\circle*{0.29}}
\multiput(10,1)(1,0){2}{\circle*{0.29}}
\put(10,0){\line(1,0){2}}
\put(10,0){\line(0,1){1}}
\put(11,0){\line(0,1){1}}
\put(11,1){\line(1,-1){1}}
\put(9.2,0.5){\makebox(0,0){$H_{3}$}}
\end{picture}\caption{Only $H_{3}$ is not a K\"{o}nig--Egerv\'{a}ry graph, as
$\alpha(H_{3})+\mu(H_{3})=4<5=\left\vert V(H_{3})\right\vert $.}%
\label{fig1}%
\end{figure}
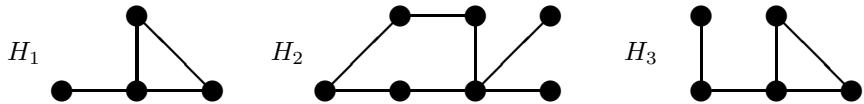

According to a well-known result of K\"{o}nig \cite{koen}, and Egerv\'{a}ry
\cite{eger}, any bipartite graph is a K\"{o}nig-Egerv\'{a}ry\emph{ }graph.
This class includes non-bipartite graphs as well (see, for instance, the
graphs $H_{1}$ and $H_{2}$ in Figure \ref{fig1}). 

It is easy to see that if $G$ is a K\"{o}nig-Egerv\'{a}ry\emph{ }graph, then
$\alpha(G)\geq\mu(G)$, and that a graph $G$ having a perfect matching is a
K\"{o}nig-Egerv\'{a}ry\emph{ }graph if and only if $\alpha(G)=\mu(G)$.

The number $d(G)=\max\{\left\vert S\right\vert -\left\vert N(S)\right\vert
:S\in\mathrm{Ind}(G)\}$ is called the \textit{critical difference} of $G$. An
independent set $A$ is \textit{critical} if $\left\vert A\right\vert
-\left\vert N(A)\right\vert =d(G)$, and the \textit{critical independence
number} $\alpha_{c}(G)$ is the cardinality of a maximum critical independent
set \cite{Zhang}. Clearly, $\alpha_{c}(G)\leq\alpha(G)$ holds for any graph
$G$. It is known that the problem of finding a critical independent set is
polynomially solvable \cite{Ageev,Zhang}.

\begin{proposition}
\label{prop4}\cite{Larson} If $S$ is a critical independent set, then there is
a matching from $N(S)$ into $S$.
\end{proposition}

If $S$ is an independent set of a graph $G$ and $H=G-S$, then we write
$G=S\ast H$. Evidently, any graph admits such representations. For instance,
if $E(H)=\emptyset$, then $G=S\ast H$ is bipartite; if $H$ is complete, then
$G=S\ast H$ is a \textit{split graph\ }\cite{FolHammer}.

\begin{proposition}
\cite{levm4}\label{prop2} $G$ is a \textit{K\"{o}nig-Egerv\'{a}ry} graph if
and only if $G=H_{1}\ast H_{2}$, where $V(H_{1})\in\Omega(G)$ and $\left\vert
V(H_{1})\right\vert \geq\mu(G)=\left\vert V(H_{2})\right\vert $.
\end{proposition}

Let $M$ be a maximum matching of a graph $G$. To adopt Edmonds's terminology
\cite{Edmonds}, we recall the following terms for $G$ relative to $M$. An
\textit{alternating path} from a vertex $x$ to a vertex $y$ is a $x,y$-path
whose edges are alternating in and not in $M$. A vertex $x$ is
\textit{exposed} relative to $M$ if $x$ is not the endpoint of a heavy edge.
An odd cycle $C$ with $V(C)=\{x_{0},x_{1},...,x_{2k}\}$ and $E(C)=\{x_{i}%
x_{i+1}:0\leq i\leq2k-1\}\cup\{x_{2k},x_{0}\}$, such that $x_{1}x_{2}%
,x_{3}x_{4},...,x_{2k-1}x_{2k}\in M$ is a \textit{blossom} relative to $M$.
The vertex $x_{0}$ is the \textit{base} of the blossom. The \textit{stem} is
an even length alternating path joining the base of a blossom and an exposed
vertex for $M$. The base is the only common vertex to the blossom and the
stem. A \textit{flower} is a blossom and its stem. A \textit{posy} consists of
two (not necessarily disjoint) blossoms joined by an odd length alternating
path whose first and last edges belong to $M$. The endpoints of the path are
exactly the bases of the two blossoms. The following result of Sterboul,
characterizes K\"{o}nig-Egerv\'{a}ry graphs in terms of forbidden configurations.

\begin{theorem}
\cite{ster} For a graph $G$, the following properties are equivalent:

\emph{(i)} $G$ is a \textit{K\"{o}nig-Egerv\'{a}ry graph};

\emph{(ii)} there exist no flower and no posy relative to some maximum
matching $M$;

\emph{(iii)} there exist no flower and no posy relative to any maximum
matching $M$.
\end{theorem}

In \cite{lov} is given a characterization of K\"{o}nig-Egerv\'{a}ry graphs
having a perfect matching, in terms of certain forbidden subgraphs with
respect to a specific perfect matching of the graph. In \cite{KoNgPeis} is
given the following characterization of K\"{o}nig-Egerv\'{a}ry graphs in terms
of excluded structures.

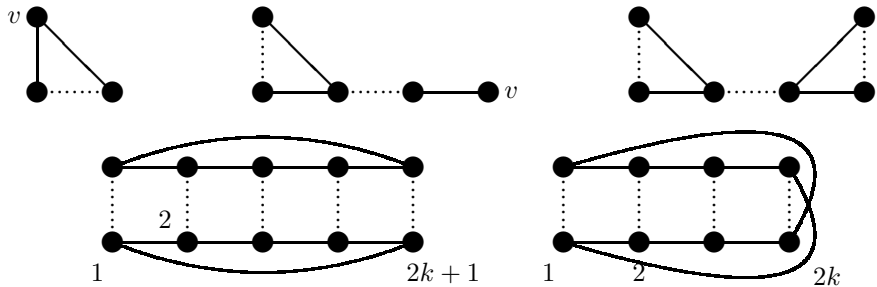
\begin{figure}[h]
\setlength{\unitlength}{1cm}\begin{picture}(5,4.3)\thicklines
\multiput(2,1)(1,0){5}{\circle*{0.29}}
\multiput(2,2)(1,0){5}{\circle*{0.29}}
\put(2,1){\line(1,0){4}}
\put(2,2){\line(1,0){4}}
\multiput(2,1)(0,0.1){10}{\circle*{0.02}}
\multiput(3,1)(0,0.1){10}{\circle*{0.02}}
\multiput(4,1)(0,0.1){10}{\circle*{0.02}}
\multiput(5,1)(0,0.1){10}{\circle*{0.02}}
\multiput(6,1)(0,0.1){10}{\circle*{0.02}}
\qbezier(2,2)(4,2.8)(6,2)
\qbezier(2,1)(4,0.2)(6,1)
\put(1.8,0.6){\makebox(0,0){$1$}}
\put(2.7,1.3){\makebox(0,0){$2$}}
\put(6.4,0.6){\makebox(0,0){$2k+1$}}
\multiput(8,2)(1,0){4}{\circle*{0.29}}
\multiput(8,1)(1,0){4}{\circle*{0.29}}
\multiput(8,1)(0,0.1){10}{\circle*{0.02}}
\multiput(9,1)(0,0.1){10}{\circle*{0.02}}
\multiput(10,1)(0,0.1){10}{\circle*{0.02}}
\multiput(11,1)(0,0.1){10}{\circle*{0.02}}
\put(8,1){\line(1,0){3}}
\put(8,2){\line(1,0){3}}
\qbezier(8,2)(12.4,3.3)(11,1)
\qbezier(8,1)(12.4,-0.3)(11,2)
\put(7.8,0.6){\makebox(0,0){$1$}}
\put(9,0.6){\makebox(0,0){$2$}}
\put(11.5,0.55){\makebox(0,0){$2k$}}
\multiput(1,3)(1,0){2}{\circle*{0.29}}
\multiput(1,3)(0.1,0){10}{\circle*{0.02}}
\put(1,4){\circle*{0.29}}
\put(1,4){\line(1,-1){1}}
\put(1,3){\line(0,1){1}}
\put(0.7,4){\makebox(0,0){$v$}}
\multiput(4,3)(1,0){4}{\circle*{0.29}}
\put(4,4){\circle*{0.29}}
\put(4,3){\line(1,0){1}}
\put(6,3){\line(1,0){1}}
\put(4,4){\line(1,-1){1}}
\multiput(4,3)(0,0.1){10}{\circle*{0.02}}
\multiput(5,3)(0.1,0){10}{\circle*{0.02}}
\put(7.3,3){\makebox(0,0){$v$}}
\multiput(9,3)(1,0){4}{\circle*{0.29}}
\multiput(9,4)(3,0){2}{\circle*{0.29}}
\multiput(9,3)(0,0.1){10}{\circle*{0.02}}
\multiput(10,3)(0.1,0){10}{\circle*{0.02}}
\multiput(12,3)(0,0.1){10}{\circle*{0.02}}
\put(9,3){\line(1,0){1}}
\put(11,3){\line(1,0){1}}
\put(9,4){\line(1,-1){1}}
\put(11,3){\line(1,1){1}}
\end{picture}\caption{Forbidden configurations. The vertex $v$ is not adjacent
to the matching edges (namely, dashed edges).}%
\label{fig44}%
\end{figure}

\begin{theorem}
\cite{KoNgPeis} Let $M$ be a maximum matching in a graph $G$. Then $G$ is a
\textit{K\"{o}nig-Egerv\'{a}ry} graph if and only if $G$ does not contain one
of the forbidden configurations, depicted in Figure \ref{fig44}, with respect
to $M$.
\end{theorem}

In \cite{Larson2009} it was shown that $G$ is a
K\"{o}nig-Egerv\'{a}ry graph if and only if $\alpha_{c}(G)=\alpha(G)$, thus
giving a positive answer to the Graffiti.pc 329 conjecture \cite{DeLaVina}.

The \textit{deficiency} of $G$, denoted by $def(G)$, is defined as the number
of exposed vertices relative to a maximum matching \cite{lovpl}. In other
words, $def(G)=\left\vert V\left(  G\right)  \right\vert -2\mu(G)$.

\mathstrut In this paper we prove that the critical difference for a
K\"{o}nig-Egerv\'{a}ry graph $G$ is given by
\[
d(G)=\left\vert \mathrm{core}(G)\right\vert -\left\vert N(\mathrm{core}%
(G))\right\vert =\alpha(G)-\mu(G)=def(G)\text{,}%
\]
and using this finding, we show that $G$ is a K\"{o}nig-Egerv\'{a}ry graph if
and only if each of its maximum independent sets is critical.\mathstrut
\mathstrut

\section{Results}

\begin{proposition}
\label{prop1}Every critical independent set is a local maximum independent set.
\end{proposition}

\begin{proof}
Suppose, on the contrary, that there is a critical independent set $S$ such
that $S\notin\Psi(G)$, i.e., there exists some independent set $A\subseteq
N[S]$, larger than $S$. It follows that $\left\vert A\cap N(S)\right\vert
>\left\vert S-S\cap A\right\vert $, and this contradicts the fact that,
according to Proposition \ref{prop4}, there is a matching from $A\cap N(S)$ to
$S$, in fact, from $A\cap N(S)$ to $S-S\cap A$.
\end{proof}

\begin{figure}[h]
\setlength{\unitlength}{1cm}\begin{picture}(5,2.2)\thicklines
\multiput(1.5,0.5)(1,0){5}{\circle*{0.29}}
\multiput(2.5,1.5)(1,0){4}{\circle*{0.29}}
\put(1.5,0.5){\line(1,0){4}}
\put(2.5,0.5){\line(1,1){1}}
\put(2.5,0.5){\line(0,1){1}}
\put(3.5,1.5){\line(1,0){1}}
\put(4.5,0.5){\line(1,1){1}}
\put(4.5,0.5){\line(0,1){1}}
\put(1.5,0.1){\makebox(0,0){$a$}}
\put(2.5,1.85){\makebox(0,0){$b$}}
\put(2.5,0.1){\makebox(0,0){$c$}}
\put(3.5,0.1){\makebox(0,0){$d$}}
\put(4.5,0.1){\makebox(0,0){$e$}}
\put(5.5,0.1){\makebox(0,0){$f$}}
\put(5.5,1.85){\makebox(0,0){$g$}}
\put(4.5,1.85){\makebox(0,0){$h$}}
\put(0.7,1){\makebox(0,0){$G_{1}$}}
\multiput(8,1)(3,0){2}{\circle*{0.29}}
\multiput(9,0)(1,0){2}{\circle*{0.29}}
\multiput(9,2)(1,0){2}{\circle*{0.29}}
\multiput(12,0)(0,1){3}{\circle*{0.29}}
\put(8,1){\line(1,0){4}}
\put(8,1){\line(1,1){1}}
\put(8,1){\line(1,-1){1}}
\put(8,1){\line(2,1){2}}
\put(8,1){\line(2,-1){2}}
\put(9,0){\line(0,1){2}}
\put(9,0){\line(1,2){1}}
\put(9,0){\line(2,1){2}}
\put(9,0){\line(1,0){1}}
\put(10,0){\line(1,1){1}}
\put(9,2){\line(1,0){1}}
\put(9,2){\line(1,-2){1}}
\put(9,2){\line(2,-1){2}}
\put(10,2){\line(1,-1){2}}
\put(10,0){\line(0,1){2}}
\put(10,0){\line(1,1){2}}
\put(11,1.3){\makebox(0,0){$v$}}
\put(12.3,0){\makebox(0,0){$x$}}
\put(12.3,1){\makebox(0,0){$y$}}
\put(12.3,2){\makebox(0,0){$z$}}
\put(7.1,1){\makebox(0,0){$G_{2}$}}
\end{picture}\caption{$\alpha(G_{1})=6$, $\mu(G_{1})=3$, $\mathrm{core}%
(G_{1})=\{a,b,d,g,f\}$ and $N(\mathrm{core}(G_{1}))=\{c,e\}$, while
$\alpha(G_{2})=4$, $\mu(G_{2})=3$, $\mathrm{core}(G_{2})=\{x,y,z\}$, and
$N(\mathrm{core}(G_{2}))=\{v\}$.}%
\label{fig112}%
\end{figure}
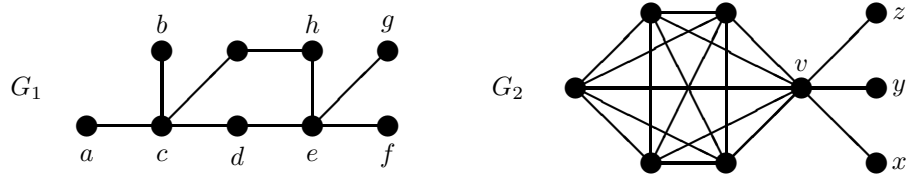

The converse of Proposition \ref{prop1} is not true; e.g., the set $\{d,h\}$
is a local maximum independent set of the graph $G_{1}$ from Figure
\ref{fig112}, but it is not critical.

Using Theorem \ref{th1}, we easily deduce the following result.

\begin{corollary}
\label{cor1}\cite{Butenko} Every critical independent set is contained in some
maximum independent set.
\end{corollary}

\begin{theorem}
\label{th2} If $G$ is a K\"{o}nig-Egerv\'{a}ry graph, then

\emph{(i)} \cite{levm4} $N(\mathrm{core}(G))=\cap\left\{  V\left(  G\right)
-S:S\in\Omega\left(  G\right)  \right\}  $;

\emph{(ii)} \cite{LevMan3} $\alpha(G)+\left\vert \cap\left\{  V\left(
G\right)  -S:S\in\Omega\left(  G\right)  \right\}  \right\vert =\mu
(G)+\left\vert \cap\left\{  S:S\in\Omega\left(  G\right)  \right\}
\right\vert $;

\emph{(iii)} \cite{LevMan3} $G-N[$\textrm{core}$(G)]$ has a perfect matching
and it is also a K\"{o}nig-Egerv\'{a}ry\emph{ }graph.
\end{theorem}

Let us notice that for non-K\"{o}nig-Egerv\'{a}ry graphs every relation
between $\alpha(G)-\mu(G)$ and $\left\vert \mathrm{core}(G)\right\vert
-\left\vert N(\mathrm{core}(G))\right\vert $ is possible.

The non-K\"{o}nig-Egerv\'{a}ry graphs from Figure \ref{fig112} satisfy:
\[
\alpha(G_{1})-\mu(G_{1})=3=\left\vert \mathrm{core}(G_{1})\right\vert
-\left\vert N(\mathrm{core}(G_{1}))\right\vert
\]
and
\[
\alpha(G_{2})-\mu(G_{2})=1<2=\left\vert \mathrm{core}(G_{2})\right\vert
-\left\vert N(\mathrm{core}(G_{2}))\right\vert .
\]

The opposite direction of the above inequality may be found in $G_{3}%
=K_{2n}-e,n\geq3$:
\[
\alpha(G_{3})-\mu(G_{3})=2-n>4-2n=2-(2n-2)=\left\vert \mathrm{core}%
(G_{3})\right\vert -\left\vert N(\mathrm{core}(G_{3}))\right\vert .
\]

\begin{theorem}
\label{th3}If $G$ is K\"{o}nig-Egerv\'{a}ry graph, then the following
equalities hold
\[
d(G)=\left\vert \mathrm{core}(G)\right\vert -\left\vert N(\mathrm{core}%
(G))\right\vert =\alpha(G)-\mu(G)=def(G)\text{.}%
\]

\end{theorem}

\begin{proof}
\textbf{\ }Firstly, let us prove that $\alpha(G)-\mu(G)\geq\left\vert
S\right\vert -\left\vert N(S)\right\vert $ holds for every $S\in
\mathrm{Ind}(G)$, i.e., $d(G)\leq\alpha(G)-\mu(G)$. If $\alpha(G)=\mu(G)$,
then $G$ has a perfect matching and
\[
\left\vert S\right\vert -\left\vert N(S)\right\vert \leq0=\alpha(G)-\mu(G)
\]
holds for every $S\in\mathrm{Ind}(G)$.

Suppose that $\alpha(G)>\mu(G)$. Let $S_{0}\in\Omega(G)$ and $M$ be a maximum
matching, i.e., $\left\vert M\right\vert =\left\vert V(G)-S_{0}\right\vert
=\mu(G)$. Assume that $S\in\mathrm{Ind}(G)$ satisfies $\left\vert S\right\vert
-\left\vert N(S)\right\vert >0$. Then one can write $S=S_{1}\cup S_{2}\cup
S_{3}$, where $S_{3}\subseteq V(G)-S_{0}$, $S_{1}\cup S_{2}\subset S_{0}$,
$S_{1}\cap S_{2}=\emptyset$, and $S_{2}$ contains every $v\in S$ matched by
$M$ with some vertex of $V(G)-S_{0}$. Since $M$ is a maximum matching, we
obtain that $\left\vert S_{2}\right\vert -\left\vert N(S_{2})\right\vert
\leq0$ and $\left\vert S_{3}\right\vert -\left\vert N(S_{3})\right\vert \leq
0$. Consequently, we infer that
\[
\alpha(G)-\mu(G)=\left\vert S_{0}\right\vert -\left\vert V(G)-S_{0}\right\vert
\geq\left\vert S_{1}\right\vert \geq\left\vert S\right\vert -\left\vert
N(S)\right\vert ,
\]
as required (see Figure \ref{fig101} for various examples of $S$%
).\begin{figure}[h]
\setlength{\unitlength}{1.0cm} \begin{picture}(5,2)\thicklines
\multiput(3,1.5)(1,0){8}{\circle*{0.29}}
\multiput(7,0.5)(1,0){4}{\circle*{0.29}}
\put(5,0.5){\circle*{0.29}}
\put(5,0.5){\line(1,0){4}}
\put(5,0.5){\line(0,1){1}}
\put(5,0.5){\line(1,1){1}}
\put(5,0.5){\line(-1,1){1}}
\put(5,0.5){\line(-2,1){2}}
\put(5,1.5){\line(2,-1){2}}
\put(7,0.5){\line(0,1){1}}
\put(7,1.5){\line(1,-1){1}}
\put(8,0.5){\line(0,1){1}}
\put(9,0.5){\line(0,1){1}}
\put(9,0.5){\line(1,1){1}}
\put(9,1.5){\line(1,-1){1}}
\put(10,0.5){\line(0,1){1}}
\put(3,1.8){\makebox(0,0){$x_{1}$}}
\put(4,1.8){\makebox(0,0){$x_{2}$}}
\put(5,1.8){\makebox(0,0){$x_{3}$}}
\put(6,1.8){\makebox(0,0){$x_{4}$}}
\put(7,1.8){\makebox(0,0){$x_{5}$}}
\put(8,1.8){\makebox(0,0){$x_{6}$}}
\put(9,1.8){\makebox(0,0){$x_{7}$}}
\put(10,1.8){\makebox(0,0){$x_{8}$}}
\put(5,0.1){\makebox(0,0){$y_{1}$}}
\put(7,0.1){\makebox(0,0){$y_{2}$}}
\put(8,0.1){\makebox(0,0){$y_{3}$}}
\put(9,0.1){\makebox(0,0){$y_{4}$}}
\put(10,0.1){\makebox(0,0){$y_{5}$}}
\end{picture}
\caption{$S_{0}=\{x_{i}:1\leq i\leq8\},M=\{y_{1}x_{4},y_{2}x_{5},y_{3}%
x_{6},y_{4}x_{7},y_{5}x_{8}\},{\ }S=S_{1}\cup S_{2}\cup S_{3}$, where
$S_{2}=\{x_{5}\},S_{3}=\{y_{4},y_{5}\}$, while $S_{1}$ belongs to
$\{\{x_{1},x_{2}\},\{x_{1}x_{3}\},\{x_{3}\}\}$.}%
\label{fig101}%
\end{figure}
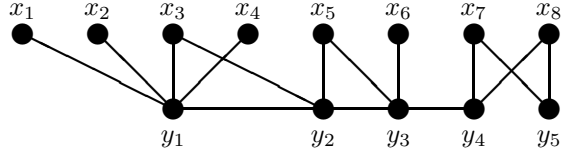

The fact that $\mathrm{core}(G)$ is an independent set of $G$ ensures that
\[
\alpha(G)-\mu(G)\geq\left\vert \mathrm{core}(G)\right\vert -\left\vert
N(\mathrm{core}(G))\right\vert .
\]

Since $G$ is a K\"{o}nig-Egerv\'{a}ry graph, we get that
\[
\alpha(G)+\mu(G)=\left\vert V\left(  G\right)  \right\vert =\left\vert
\mathrm{core}(G)\right\vert +\left\vert N(\mathrm{core}(G))\right\vert
+\left\vert V(G-N[\mathrm{core}(G)])\right\vert .
\]

Assuming that
\[
\alpha(G)-\mu(G)>\left\vert \mathrm{core}(G)\right\vert -\left\vert
N(\mathrm{core}(G))\right\vert ,
\]
we obtain the following contradiction%
\begin{align*}
2\alpha(G)  &  >2\left\vert \mathrm{core}(G)\right\vert +\left\vert
V(G-N[\mathrm{core}(G)])\right\vert \\
&  =2\left\vert \mathrm{core}(G)\right\vert +2\alpha(G-N[\mathrm{core}%
(G)])=2\alpha(G),
\end{align*}
because $\left\vert V(G-N[\mathrm{core}(G)])\right\vert =2\alpha\left(
G-N[\mathrm{core}(G)]\right)  $ by Theorem \ref{th2}\emph{(iii)}.

Therefore, we get that $\alpha(G)-\mu(G)=\left\vert \mathrm{core}%
(G)\right\vert -\left\vert N(\mathrm{core}(G))\right\vert $. Actually, this
equality immediately follows from Theorem \ref{th2}\emph{(i)},\emph{(ii)}, but
the current way of proof exploits different aspects of $\mathrm{Ind}(G)$.

Further, using the inequality $d(G)\leq\alpha(G)-\mu(G)$ and the equality
\[
\alpha(G)-\mu(G)=\left\vert \mathrm{core}(G)\right\vert -\left\vert
N(\mathrm{core}(G))\right\vert ,
\]
we finally deduce that%
\begin{align*}
\left\vert \mathrm{core}(G)\right\vert -\left\vert N(\mathrm{core}%
(G))\right\vert  &  \leq\max\{\left\vert S\right\vert -\left\vert
N(S)\right\vert :S\in\mathrm{Ind}(G)\}=d(G)\\
&  \leq\alpha(G)-\mu(G)=\left\vert \mathrm{core}(G)\right\vert -\left\vert
N(\mathrm{core}(G))\right\vert ,
\end{align*}
i.e.,
\[
\alpha(G)-\mu(G)=\left\vert \mathrm{core}(G)\right\vert -\left\vert
N(\mathrm{core}(G))\right\vert =d\left(  G\right)  .
\]

Since $G$ is a K\"{o}nig-Egerv\'{a}ry graph, we infer that
\[
\alpha(G)-\mu(G)=\alpha(G)+\mu(G)-2\mu(G)=\left\vert V\left(  G\right)
\right\vert -2\mu(G)=def(G),
\]
and this completes the proof.
\end{proof}

\begin{corollary}
If $G$ is a K\"{o}nig-Egerv\'{a}ry graph, then $d(G)=0$ if and only if $G$ has
a perfect matching.
\end{corollary}

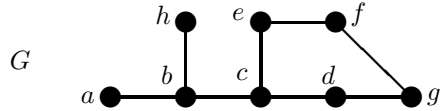
\begin{figure}[h]
\setlength{\unitlength}{1.0cm} \begin{picture}(5,1.2)\thicklines
\multiput(4.5,0)(1,0){5}{\circle*{0.29}}
\multiput(5.5,1)(1,0){3}{\circle*{0.29}}
\put(4.5,0){\line(1,0){4}}
\put(5.5,0){\line(0,1){1}}
\put(6.5,1){\line(1,0){1}}
\put(6.5,0){\line(0,1){1}}
\put(7.5,1){\line(1,-1){1}}
\put(4.2,0){\makebox(0,0){$a$}}
\put(5.25,0.3){\makebox(0,0){$b$}}
\put(5.2,1.1){\makebox(0,0){$h$}}
\put(6.25,0.3){\makebox(0,0){$c$}}
\put(7.4,0.3){\makebox(0,0){$d$}}
\put(6.2,1.1){\makebox(0,0){$e$}}
\put(7.8,1.1){\makebox(0,0){$f$}}
\put(8.8,0){\makebox(0,0){$g$}}
\put(3.3,0.5){\makebox(0,0){$G$}}
\end{picture}\caption{$G$ has $\alpha(G)=4,\mu(G)=3$,{ }$\mathrm{core}%
(G)=\{a,h\}$ and $N\left(  \mathrm{core}\left(  G\right)  \right)  {=\{b\}}$.}%
\label{fig10}%
\end{figure}

\begin{remark}
There exist non-K\"{o}nig-Egerv\'{a}ry graphs enjoying the equalities
\[
d(G)=\left\vert \mathrm{core}(G)\right\vert -\left\vert N(\mathrm{core}%
(G))\right\vert =\alpha(G)-\mu(G),
\]
see, for instance, the graph $G$ from Figure \ref{fig10}.
\end{remark}

\begin{theorem}
The following assertions are equivalent:

\emph{(i)} $G$ is a K\"{o}nig-Egerv\'{a}ry graph;

\emph{(ii)} there is $S\in\Omega(G)$, such that $S$ is critical, i.e.,
$\alpha_{c}(G)=\alpha(G)$;

\emph{(iii)} every $S\in\Omega(G)$ is critical.
\end{theorem}

\begin{proof}
\emph{(i)} $\Longrightarrow$\ \emph{(iii) }Let $S\in\Omega(G)$, $A=S-$%
\textrm{core}$(G)$ and $B=V\left(  G\right)  -S-N(\mathrm{core}(G))$. By
Theorem\emph{ }\ref{th2}\emph{(iii)}, we infer that $\left\vert A\right\vert
=\left\vert B\right\vert $, since $G-N[$\textrm{core}$(G)]$ has a perfect
matching. Hence, we obtain that
\begin{align*}
\left\vert S\right\vert -\left\vert N(S)\right\vert  &  =\left\vert
A\right\vert +\left\vert \mathrm{core}(G)\right\vert -(\left\vert B\right\vert
+\left\vert N(\mathrm{core}(G))\right\vert \\
&  =\left\vert \mathrm{core}(G)\right\vert -\left\vert N(\mathrm{core}%
(G))\right\vert .
\end{align*}
In other words, according to Theorem \ref{th3}, the equality $\left\vert
S\right\vert -\left\vert N(S)\right\vert =d(G)$ is true for every $S\in
\Omega(G)$.

\emph{(iii)} $\Longrightarrow$\ \emph{(ii)} It is clear.

\emph{(ii)} $\Longrightarrow$\ \emph{(i)} This was done in\emph{
}\cite{Larson2009}. For the sake of completeness we add the proof.

There is a critical independent set $S$ with $\left\vert S\right\vert
=\alpha_{c}(G)=\alpha(G)$. By Proposition \ref{prop4}, there exists a matching
$M$ from $N(S)$ into $S$, and clearly, $\left\vert M\right\vert =\left\vert
N(S)\right\vert =\mu(G)$. Hence, we finally obtain that $\left\vert
V(G)\right\vert =\left\vert S\right\vert +\left\vert N(S)\right\vert
=\alpha(G)+\mu(G)$, i.e., $G$ is a K\"{o}nig-Egerv\'{a}ry graph.
\end{proof}

\section{Conclusions}

In this paper we give a new characterization of K\"{o}nig-Egerv\'{a}ry\emph{
}graphs. On the one hand, it is similar in form to Sterboul's theorem
\cite{ster}. On the other hand it extends Larson's finding \cite{Larson2009}.
We found that the critical difference of a K\"{o}nig-Egerv\'{a}ry\emph{ }graph
$G$ is given by
\[
d(G)=\left\vert \mathrm{core}(G)\right\vert -\left\vert N(\mathrm{core}%
(G))\right\vert =\alpha(G)-\mu(G)=def(G).
\]
It seems interesting to find other families of graphs satisfying these equalities.


\begin{thebibliography}{99}                                                                                               %


\bibitem {Ageev}A. A. Ageev, \emph{On finding critical independent and vertex
sets}, SIAM J. Discrete Mathematics \textbf{7} (1994) 293--295.

\bibitem {bourhams1}J. - M. Bourjolly, P. L. Hammer, B. Simeone, \emph{Node
weighted graphs having K\"{o}nig-Egervary property}, Math. Programming Study
\textbf{22} (1984) 44-63.

\bibitem {bourpull}J. M. Bourjolly, W. R. Pulleyblank,
\emph{K\"{o}nig-Egerv\'{a}ry graphs, 2-bicritical graphs and fractional
matchings}, Discrete Applied Mathematics \textbf{24} (1989) 63--82.

\bibitem {Butenko}S. Butenko, S. Trukhanov, \emph{Using Critical Sets to Solve
the Maximum Independent Set Problem}, Operations Research Letters \textbf{35}
(2007) 519-524.

\bibitem {DeLaVina}E. DeLaVina, \emph{Written on the Wall II, Conjectures of
Graffiti.pc},\newline http://cms.dt.uh.edu/faculty/delavinae/research/wowII/

\bibitem {dem}R. W. Deming, \emph{Independence numbers of graphs - an
extension of the K\"{o}nig-Egerv\'{a}ry theorem}, Discrete Mathematics
\textbf{27} (1979) 23--33.

\bibitem {Edmonds}J. Edmonds, \emph{Paths, trees and flowers}, Canadian
Journal of Mathematics \textbf{17} (1965) 449-467.

\bibitem {eger}E. Egerv\'{a}ry, \emph{On combinatorial properties of
matrices}, Matematikai Lapok \textbf{38} (1931) 16--28.

\bibitem {FolHammer}S. F\"{o}ldes, P. L. Hammer, \emph{Split graphs},
Proceedings of 8th Southeastern Conference on Combinatorics, Graph Theory and
Computing (F. Hoffman \textit{et al.} eds), Louisiana State University, Baton
Rouge, Louisiana, 311--315.

\bibitem {koen}D. K\"{o}nig, \emph{Graphen und Matrizen}, Matematikai Lapok
\textbf{38} (1931) 116--119.

\bibitem {korach}E. Korach, \emph{On dual integrality, \textit{min-max}%
}\textit{\ }\emph{equalities and algorithms in combinatorial programming},
University of Waterloo, Department of Combinatorics and Optimization, Ph.D.
Thesis, 1982.

\bibitem {KoNgPeis}E. Korach, T. Nguyen, B. Peis, \emph{Subgraph
characterization of Red/Blue-Split Graph and K\"{o}nig-Egerv\'{a}ry\ graphs},
Proceedings of the Seventeenth Annual ACM-SIAM Symposium on Discrete
Algorithms, ACM Press (2006) 842-850.

\bibitem {Larson}C. E. Larson, \emph{A note on critical independence
reductions}, Bulletin of the Institute of Combinatorics and its Applications 5
(2007) 34-46.

\bibitem {Larson2009}C. E. Larson, \emph{A new characterization of
K\"{o}nig-Egerv\'{a}ry graphs}, The 2$^{nd}$ Canadian Discrete and Algorithmic
Mathematics Conference, May 25-28, 2009, CRM Montreal (Canada).

\bibitem {levm2}V. E. Levit, E. Mandrescu, \emph{Well-covered and
K\"{o}nig-Egerv\'{a}ry graphs}, Congressus Numerantium \textbf{130} (1998) 209--218.

\bibitem {LevMan2}V. E. Levit, E. Mandrescu, \emph{A new greedoid: the family
of local maximum stable sets of a forest}, Discrete Applied Mathematics
\textbf{124} (2002) 91-101.

\bibitem {levm3}V. E. Levit, E. Mandrescu, \emph{Combinatorial properties of
the family of maximum stable sets of a graph}, Discrete Applied Mathematics
\textbf{117} (2002) 149-161.

\bibitem {levm4}V. E. Levit, E. Mandrescu, \emph{On }$\alpha^{+}$\emph{-stable
K\"{o}nig-Egerv\'{a}ry graphs}, Discrete Mathematics \textbf{263} (2003) 179--190.

\bibitem {LevMan3}V. E. Levit, E. Mandrescu, \emph{On }$\alpha$\emph{-critical
edges in K\"{o}nig-Egerv\'{a}ry graphs},\emph{\ }Discrete Mathematics
\textbf{306} (2006) 1684-1693.

\bibitem {lov}L. Lov\'{a}sz, \emph{Ear decomposition of matching covered
graphs}, \textit{Combinatorica }\textbf{3} (1983) \ 105-117.

\bibitem {lovpl}L. Lov\'{a}sz, M. D. Plummer, \emph{Matching Theory}, Annals
of Discrete Mathematics \textbf{29} (1986) North-Holland.

\bibitem {NemhTro}G. L. Nemhauser and L. E. Trotter, Jr., \emph{Vertex
packings: structural properties and algorithms}, Mathematical Programming
\textbf{8} (1975) 232-248.

\bibitem {pasdema}V. T. Paschos, M. Demange, \emph{A generalization of
K\"{o}nig-Egerv\'{a}ry graphs and heuristics for the maximum independent set
problem with improved approximation ratios}, European Journal of Operational
Research \textbf{97} (1997) 580--592.

\bibitem {pulleybl}W. R. Pulleyblank, \emph{Matchings and Extensions}, in:
\textit{Handbook of Combinatorics, Volume 1} (eds. R. L. Graham, M. Grotschel
and L. Lov\'{a}sz), MIT Press and North-Holland, Amsterdam (1995) 179-232.

\bibitem {ster}F. Sterboul, \emph{A characterization of the graphs in which
the transversal number equals the matching number}, Journal of Combinatorial
Theory Series B \textbf{27} (1979) 228--229.

\bibitem {Zhang}C. Q. Zhang, \emph{Finding critical independent sets and
critical vertex subsets are polynomial problems}, SIAM J. Discrete Mathematics
\textbf{3} (1990) 431-438.
\end{thebibliography}
\end{document}